\documentclass{article}

\usepackage{amsmath,amsthm,amssymb,amscd}
\usepackage{amsfonts}
\usepackage{rotating}
\usepackage{euscript}
\usepackage{pst-node}
\usepackage{epsfig,verbatim}
\usepackage{pictexwd,dcpic}

\def\be{\begin{equation}}
\def\ee{\end{equation}}

\def\C{{\mathbb C}} 
\def\f{\EuScript}
 
\def\P{{\mathbb P}}

\def\ord{{\rm ord}}

\def\phi{{\varphi}}
\def\v{{\varepsilon}} 
\def\tt{\widetilde}
\def\deg{{\rm deg\,}}

\def\cos{{\rm cos\,}}

\def\GCD{{\rm GCD }}
\def\LCM{{\rm LCM }}

\def\bp{\begin{proposition}}
\def\ep{\end{proposition}}
\def\rad{{\mathrm{rad}}}

\def\bt{\begin{theorem}}
\def\et{\end{theorem}}
\def\br{\begin{remark}}
\def\er{\end{remark}}
\def\be{\begin{equation}}
\def\bee{\begin{equation*}}
\def\l{\label}
\def\la{\label}

\def\ee{\end{equation}}
\def\eee{\end{equation*}}
\def\bl{\begin{lemma}}
\def\el{\end{lemma}}
\def\bc{\begin{corollary}}
\def\ec{\end{corollary}}
\def\pr{\noindent{\it Proof. }}

\def\bd{\begin{definition}}
\def\ed{\end{definition}}
\def\t{\widetilde}
\def\tilde{\widetilde}

\newtheorem{theorem}{Theorem}[section]
\newtheorem{lemma}[theorem]{Lemma}
\newtheorem{corollary}[theorem]{Corollary}
\newtheorem{proposition}[theorem]{Proposition}
\theoremstyle{definition}
\newtheorem{remark}[theorem]{Remark}

\begin{document}
\title{Polynomial semiconjugacies, decompositions of iterations, and invariant curves }
\author{F. Pakovich}
\maketitle
\begin{abstract}
We study 
the functional equation $A\circ X=X\circ B$, 
where $A,$ $B$, and $X$ 
are polynomials with complex coefficients. 
Using results of \cite{p1} about polynomials sharing preimages of compact sets in $\C$, we show that 
for given $B$ its solu\-tions may be described in terms 
of the filled-in Julia set of  $B$. 
On this base, we prove a number of results describing    
a general structure of solutions.
The results obtained imply in particular the result of Medvedev and 
Scanlon \cite{ms} about invariant curves of maps $F:\,\C^2 \rightarrow \C^2$ 
of the form $(x,y)\rightarrow (f(x),f(y))$, where $f$ is a polynomial,
and a version of the result of Zieve and M\"uller \cite{mz} about decompositions 
of iterations of a polynomial.

\end{abstract}

\begin{section}{Introduction}
Let $A$ and $B$ be rational functions of degree at least two on the Riemann sphere. The functions $A$ and $B$ are called commuting if 
\be \l{0} A\circ B=B\circ A,\ee
and conjugate if 
\be \l{1}
A\circ X=X\circ B
\ee for some rational function $X$ of degree one. 

In case if \eqref{1} is satisfied for some rational function $X$ of degree at least two, 
the function $B$ is called semiconjugate to $A$, and the function $X$ is called a semiconjugacy from 
$B$ to $A.$ 
In distinction with the conjugation, the semiconjugation is not an equivalency relation.  We will use the notation 
$A\leq B$
if for given rational functions $A$ and $B$ there exists a non-constant rational function $X$ such that \eqref{1} holds, and the notation  
$A\underset{X}{\leq} B$  if $A$,$B$, and $X$ satisfy 
\eqref{1}.
The notation reflects the fact that the binary relation on the set of rational functions 
defined by equality \eqref{1} 
 is a  preorder. Indeed, it follows from  
$A\underset{X}{\leq} B$  and
$B\underset{Y}{\leq} C$ that 
$A\underset{X\circ Y}{\leq} C$.

Both equations \eqref{0} and \eqref{1} have ``obvious'' solutions. Namely, equation \eqref{0} has solutions  of the form 
\be \la{rel1} A=R^{\circ m}, \ \ \ \ B=R^{\circ n},
\ee
where $R$ is an arbitrary rational function and $m,n\geq 1.$ Notice that 
such $A$ and $B$ have an iteration in common, that is
\be \la{rel2} A^{\circ n}=B^{\circ m}
\ee
for some $n,m\geq 1.$

In order to obtain solutions of equation \eqref{1} we can take arbitrary rational functions $A_1,B_1$ and set 
$$F=A_1\circ B_1, \ \ \ G=B_1\circ A_1.$$ 
Then the equality 
\be \l{ii} (A_1\circ B_1)\circ A_1= A_1\circ (B_1\circ A_1)\ee
implies that $F\underset{A_1}{\leq} G$. Similarly, 
$G\underset{B_1}{\leq} F$. Moreover, if  now
$A_2,B_2$ are rational functions such that the equality 
\be \l{ko}  G=A_{2}\circ B_{2}\ee holds, 
then the function $H=B_2\circ A_2$ satisfies $G\underset{A_2}{\leq} H$ and
$H\underset{B_2}{\leq} G$, implying that $F\underset{A_1\circ A_2}{\leq} H$ and
$H\underset{B_2\circ B_1}{\leq} F.$ 
This motivates the following definition of an equivalency relation on the set of rational functions:
$F\sim G$ if there exist rational functions $A_i,B_i,$ $1\leq i \leq n,$ such that $$F=A_1\circ B_1, \ \ \ \ \  G=B_n\circ A_n,$$ and 
$$B_i\circ A_i=A_{i+1}\circ B_{i+1}, \ \ \ 1\leq i \leq n-1.$$ 
 Clearly, $F\sim G$ implies that $F\leq G$ and $G\leq F$.  
Notice that since for any rational function $X$ of degree one the equality
$$A=(A\circ X)\circ X^{-1}$$ implies that $A\sim X^{-1}\circ A\circ X$, any equivalence class  is a collection of conjugacy classes.

Functional equation \eqref{0} was first studied by Fatou, Julia, and Ritt in the papers 
\cite{f}, \cite{j}, and \cite{r}. 
In all these papers it was  
assumed that the considered commuting functions $A$ and $B$ have no 
iterate in common.
Fatou and Julia described solutions of \eqref{0} under the additional  
assumption that the
Julia set of $A$ or $B$ does not coincide with the whole complex plane, and Ritt investigated the general case.
Briefly, the Ritt theorem states that if rational functions $A$ and $B$ commute and
no iterate of $A$ is equal to an iterate of $B$,
then, up to a conjugacy, $A$ and $B$ are either powers, or  Chebyshev polynomials, or Latt\`es functions. 
Another proof of the Ritt theorem
was given by Eremenko in \cite{e2}. 
Notice however that a description of commuting $A$ and $B$ with a common iterate
is known only in the polynomial case.  Thus, in a certain sense the classification of commuting rational functions is not yet completed. 
On the other hand, it was shown by Ritt  (\cite{r3}, \cite{r}) that in the polynomial case
equality \eqref{0} implies that, up to the change
$$ A\rightarrow \lambda \circ A\circ  \lambda^{-1},\ \ \ B\rightarrow \lambda \circ B\circ  \lambda^{-1},$$ where $\lambda$ is a polynomial of degree one, either 
$$ A=z^n, \ \ \ B=\v z^m,$$ where $\v^n=\v,$ or
$$ A=\pm T_n, \ \ \ B=\pm T_m,$$
or 
$$ A= \v_1R^{\circ m}, \ \ \ \ B=\v_2R^{\circ n},$$
where $R=zS(z^{\ell})$ for some polynomial $S$ and  $\v_1,$ $\v_2$ are $l$-th roots of unity. 
In fact, this conclusion remains true if instead of \eqref{0} one were to assume only that  $A$ and $B$ 
share a completely invariant compact set in $\C$ (see \cite{p1}).

Equation \eqref{1} was investigated in the recent paper \cite{paksem}.
The main result of \cite{paksem} states that if a rational function $B$ is semiconjugate to a rational function $A$, then either $A\sim B$, or 
$A$ and $B$ are 
``minimal holomorphic self-maps'' between orbifolds of non-negative Euler characteristic on the Riemann sphere. 
The last class of functions is a natural  extension of the class of Latt\`es functions and admits a neat characterization.
However, similar to the description of commuting rational functions, the description of solutions of 
\eqref{1} given in \cite{paksem}
is not completely satisfactory, since it gives no information about  equivalent rational functions. 
In particular, the results of \cite{paksem} do not provide any bounds on the number of conjugacy classes in an equivalence class of a rational function $B$ or more generally on the number of conjugacy classes of $A$ such that $A\leq B$. 
Another related problem is the following: 
is it  true that if conditions $A\leq B$ and $B\leq A$ hold simultaneously, then 
$A\sim B$ ?  
Finally,
it would be desirable to obtain some handy structural descriptions of the totality of $X$ satisfying \eqref{1} for given $A$ and $B$, and of
the totality of $A$ satisfying $A\leq B$ for given $B$.

In this paper we study equation \eqref{1} with emphasis on the above questions in the case where all the functions involved are {\it polynomials}. 
Notice that in distinction with the general case, for  polynomials there exists quite a comprehensive theory of functional decompositions developed by Ritt \cite{r1}. Nevertheless, questions regarding polynomial decompositions 
may be highly non-trivial, and a number of recent papers are devoted to such questions arising from different branches of mathematics. 
Let us mention for example the paper \cite{mz} with applications to algebraic dynamics (\cite{gtz2}), or the paper  
\cite{pakmom} with applications to differential equations (\cite{paksol}). Another example is the recent paper \cite{ms} about invariant varieties for dynamical systems defined by coordinatwise actions of polynomials,
a considerable part of which
concerns properties of polynomial solutions of \eqref{1}.  

The main distinction between this paper and the above mentioned papers 
is the systematical use of ideas and results from the paper \cite{p1} which relates 
polynomials sharing preimages of 
compact sets in $\C$ with the functional equation 
$A\circ C=D\circ B$. In particular, the main result of \cite{p1} leads to a characterization of polynomial solutions of \eqref{1} in terms of filled-in Julia sets.
Recall that for a polynomial $B$ the filled-in Julia set  $K(B)$ is defined as the set of points in $\C$ whose orbits under iterations of $B$ are bounded.
Since equality \eqref{1} implies the equalities  
$$A^{\circ n}\circ X=X\circ B^{\circ n}, \ \ \  n\geq 1,$$ 
it it easy to see that if $X$ is a semiconjugacy from $B$ to $A$, then the preimage 
$X^{-1}(K(A))$ coincides with $K(B)$.
We show that this property is in fact characteristic.

\bt \l{j+} Let $A$, $B$ and $X$ be polynomials of degree at least two such that $A\underset{X}{\leq} B$.
Then
 \be \l{osn} X^{-1}(K(A))=K(B).\ee
In the other direction, if equality \eqref{osn} holds and $\deg A=\deg B$, then there exists a polynomial of degree one 
$\mu$ such that 
$$(\mu \circ A)\circ X=X\circ B$$ and 
$\mu(K(A))=K(A).$ More generally, if for given $B$ and $X$  the condition \be \l{blia} X^{-1}(K)=K(B)\ee
holds for some compact set $K$ in $\C$, then there exists a polynomial $A$ such that 
$A\underset{X}{\leq} B$
and $K(A)=K$.
\et

For a fixed polynomial $B$ of degree at least two denote by $\f E(B)$ the set of  polynomials $X$ of degree at least two such that 
$A\underset{X}{\leq} B$ for some polynomial $A$. An immediate corollary of Theorem \ref{j+} is that a polynomial $X$ is contained in $\f E(B)$ if and only if $K(B)$ is a union of fibers of $X$. 
Another corollary is that if $A\underset{X}{\leq} B$, then for any decomposition  $X=X_1\circ X_2$  there exists 
a polynomial $C$ such that 
$$A\underset{X_1}{\leq} C, \ \ \ \ C\underset{X_2}{\leq} B.$$
Notice that in particular this puts  
the problem of the description of decompositions of iterations of a polynomial,   first considered in the paper \cite{mz}, into the context of equation \eqref{1}.
Indeed, since $B\circ B^{\circ d}=B^{\circ d}\circ B,$ the polynomial $B^{\circ d}$ is contained in $\f E(B)$ and hence 
for any decomposition 
$B^{\circ d}=Y\circ X$   
the equalities $$B\circ Y=Y\circ A, \ \ \ \ A\circ X=X\circ B$$ hold for some polynomial $A$.

\vskip 0.1cm
The following statement also is a corollary of the main result of \cite{p1}.

\bt \l{resh} For any $X_1$, $X_2\in \f E(B)$ there exists
$X\in \f E(B)$ such that  $\deg X=\LCM(\deg X_1,\deg X_2)$ and $$ X=U_1\circ X_1=U_2\circ X_2$$ for some polynomials $U_1,$ $U_2$. Furthermore,  
there exists  
$W\in \f E(B)$ such that 
$\deg W=\GCD(\deg X_1,\deg X_2)$ and 
$$ X_1=V_1\circ W, \ \ \ X_2=V_2\circ W$$ for some polynomials $V_1,$ $V_2$. 
\et

For fixed polynomials $A$, $B$ denote  by $\f E(A,B)$ the subset  of $\f E(B)$  (possibly empty) consisting of polynomials $X$ such that 
$A\underset{X}{\leq} B$.  In particular, the set  $\f E(B,B)$ consists of polynomials of degree at least two commuting with $B$.
We will call a polynomial $P$ {\it special} if it is conjugated to $z^n$ or $\pm T_n$, or
equivalently if there exists a M\"obius transformation $\mu$ which 
maps $K(P)$ to $\mathbb D$ or $[-1,1].$ 
The following result describes a general structure  of $\f E(A,B)$ for non-special 
$A$, $B.$

\bt \l{uni} Let $A$ and $B$ be fixed  non-special polynomials of degree at least two
such that the set $\f E(A,B)$ is non-empty, and let $X_0$ be an element of  $\f E(A,B)$ of the minimum possible degree. Then a polynomial $X$ belongs to 
$\f E(A,B)$ if and only if $X= \t A\circ X_0$ for some polynomial  $\t A$ commuting with $A.$
\et

Notice that in a sense this result is a generalization of the result of Ritt about commuting polynomials. Indeed, applying Theorem \ref{uni} for $B=A$ and $X=B$, we obtain that if
 $A$ is non-special and $B\in \f E(A,A)$, then $B=\t A\circ R$, where $R$ is a polynomial of the minimum possible degree in $\f E(A,A)$.
Now we can apply Theorem \ref{uni} again to the polynomial $\t A$ and so on, arriving eventually to the representation 
$B=\mu_1 \circ R^{\circ m_1}$, where $\mu_1$ is a polynomial  of degree one commuting with $A$. 
In particular, since $A\in \f E(A,A)$, the equality  
$A=\mu_2\circ R^{\circ m_2}$ holds 
for some polynomial  $\mu_2$ of degree one  commuting with $A$.

Another corollary of Theorem \ref{uni} is the following result obtained by Medve\-dev and Scanlon in  the paper \cite{ms}: if $\f C\subset \C^2$ is an irreducible algebraic curve invariant  under the map $F:\,(x,y)\rightarrow (f(x),f(y)),$ where $f$ is a non-special polynomial, then there exists a polynomial $p$ which commutes with $f$ such that 
$\f C$ has the form $z_1= p(z_2)$ or $z_2=p(z_1)$. More general, we prove the following statement which supplements the results of \cite{ms} about algebraic curves invariant  under the map $F:\,(x,y)\rightarrow (f(x),g(y)),$ where $f$ and $g$ are non-special polynomials.

\bt \l{new} Let $f$ and $g$ be non-special polynomials of degree at least two and $\f C$ a curve in $\C^2$. Then $\f C$ is an irreducible $(f,g)$-invariant curve if and only if 
 $\f C$ has the form $u(x)-v(y)=0$, where $u,v$  are  polynomials of coprime degrees satisfying a system given by the equations 
\be \l{krys} t\circ u=u\circ f, \ \ \ \ t\circ v=v\circ g\ee for some polynomial $t.$
\et

\vskip 0.2cm

Our next result 
describes the interrelations between the equivalence $\sim$, the preorder $\leq\,$, and  
decompositions of iterations.

\bt \l{e}
Let $A$ and $B$ be polynomials of degree at least two. Then conditions $A\leq B$ and $B\leq A$  hold simultaneously
if and only if  $A\sim B$.
Furthermore, 
$A\sim B$ if and only if there exist 
polynomials $X$, $Y$ such that 
$$B\circ Y=Y\circ A, \ \ \ \ A\circ X=X\circ B,$$
and $Y\circ X=B^{\circ d}$ for some $d\geq 0.$
\et

For a fixed polynomial $B$ of degree at least two  denote by $\f F(B)$ the set  of polynomials $A$ such that 
$A\leq B$. 
The following theorem gives a structural description of  the set $\f F(B)$.

\bt \l{exi}
Let $B$ be a fixed  non-special polynomial of degree $n\geq 2.$ 
Then there exist $A\in \f F(B)$ and a semiconjugacy $X$ from $B$ to $A$ which are universal in the following sense:
for any polynomial $C\in \f F(B)$ 
there exist polynomials $X_C,$ $U_C$ such that
$X=U_C\circ X_C$ and
the diagram  
\be 
\begin{CD} \l{gooopa2}
\C @>B>> \C\\
@VV X_C V @VV X_C V\\ 
\C @>C>> \C\\
@VV U_C V @VV U_C V\\ 
\C @>  A >> \C\ 
\end{CD}
\ee
is commutative. Furthermore, 
the degree of $X$ is 
bounded from above by a constant $c=c(n)$ which depends  
on $n$ only.
\et 

We did not make special efforts to obtain an optimal estimation for $c(n),$ however our method of proof
shows that $$c(n)\leq (n-1)!n^{2\log_2n+3}.$$  
Thus, Theorem \ref{exi}  gives an effective bound on the number of conjugacy classes of polynomials $A$ such that $A\leq B$. 


\vskip 0.2cm

The paper is organized as follows. In the second section we give a very brief overview of the Ritt theory. 
In the third section we recall basic results of \cite{p1} and prove Theorem \ref{j+} and Theorem \ref{resh}. We also prove the corollaries of Theorem \ref{j+} mentioned above.  
In the fourth section we first show that 
if $A\leq B$ and one of polynomials $A$ or $B$ is special, then the other one also is  
special (Theorem \ref{sem}). Then we prove Theorem \ref{uni} and deduce from it the result of Ritt about commuting polynomials. We also 
apply Theorem \ref{uni} to the problem  of description of curves in $\C^2$  
invariant under maps $F:\,(x,y)\rightarrow (f(x),g(y)),$ where $f$ and $g$ are polynomials, and prove Theorem \ref{new}. Finally, we
prove Theorem \ref{e}.


In the fifth section 
we first show (Theorem \ref{ui}) that if $B$ is a non-special polynomial of degree $n$, and $X\in \f E(B)$, then the degree $l$ of any special compositional factor of $X$ satisfies the inequality $l\leq 2n.$  
On this base  we prove that if  $X\in\f E(B)$ is not a polynomial in $B$, 
then $\deg X$ is bounded from above by a constant  which depends  
on $n$ only. 
In turn, from this result we deduce Theorem \ref{exi}. 
As another corollary of the boundedness of $\deg X$ we obtain the 
following  result of  Zieve and M\"uller (\cite{mz}): 
if  $B$ is a non-special polynomial of degree $n\geq 2$, 
and $X$ and $Y$ are polynomials such that $Y\circ X=B^{\circ s}$ for some $s\geq 1$, then there exist polynomials $\t X$, $\t Y$ and $i,j\geq 0$ such that 
$$Y=B^{\circ i}\circ \t Y, \ \ \ X=\t X\circ B^{\circ j},\ \ \ {\rm and} \ \ \ \t Y\circ \t X=B^{\circ \t s},$$   
where $\t s$ is bounded from above by a constant  which depends
on $n$ only. 
\end{section}

\begin{section}{\l{overv} Overview of the Ritt theory}
Let $F$ be a polynomial with complex coefficients.
The polynomial $F$ is called {\it indecomposable} if
the equality $F=F_2\circ F_1$ implies that at least one of the polynomials $F_1,F_2$ is of degree one. 
Any representation 
of a polynomial $F$ in the form $F=F_r\circ F_{r-1}\circ \dots \circ F_1,$
where $F_1,F_2,\dots, F_r$ are polynomials,
is called {\it a decomposition} of $F.$ A decomposition is called {\it maximal}
if all $F_1,F_2,\dots, F_r$ are
indecomposable and of degree greater than one.
Two decompositions having an equal number of terms
$$F=F_r\circ F_{r-1}\circ \dots \circ F_1 \ \ \ \ {\rm and} \ \ \ \  
F=G_{r}\circ G_{r-1}\circ \dots \circ G_1$$
are called {\it equivalent} if either $r=1$ and $F_1=G_1$, or $r\geq 2$ and there exist polynomials  $\mu_i,$ $1\leq i \leq r-1,$ of degree 1 such that 
$$F_r=G_r\circ \mu_{r-1}, \ \ \ 
F_i=\mu_{i}^{-1}\circ G_i \circ \mu_{i-1}, \ \ \ 1<i< r, \ \ \ {\rm and} \ \ \ F_1=\mu_{1}^{-1}\circ G_1.
$$

The theory of polynomial decompositions established by  Ritt can be summarized in the form
of two theorems usually called the first and the second Ritt
theorems (see \cite{r1}). 

The first Ritt theorem states roughly speaking that any  maximal decompositions of a polynomial may be obtained from any other by some iterative process involving the functional equation 
\be \la{ura} A\circ C=D\circ B.\ee 

\bt [\cite{r1}]\la{rr1} Any two maximal decompositions $\f D, \f E$ 
of a 
polynomial $P$ have an equal number of terms. Furthermore, there exists 
a chain of maximal decompositions $\f F_i$, $1\leq i \leq s,$ of $P$ such that 
$\f F_1=\f D,$ $\f F_s\sim \f E,$ and $\f F_{i+1}$ is obtained from $\f F_i$ 
by a replacement of  
two successive polynomials $A\circ C$ in $\f F_i$
by two other polynomials $D\circ B$ 
such that \eqref{ura} holds.
\et

The second Ritt theorem in turn describes indecomposable  polynomial solutions of \eqref{ura}. More precisely,  
it describes solutions satisfying the condition
\be \l{ggccdd} {\rm GCD}(\deg A, \deg D)=1, \ \ \ {\rm GCD}(\deg C,\deg B)=1,\ee
which holds in particular if $A,C,D,B$ are indecomposable (see Theorem \ref{r1} below).

\bt [\cite{r1}]\la{r2}
Let $A,C,D,B$ be polynomials  
such that \eqref{ura} and \eqref{ggccdd} hold.
Then there exist polynomials $\sigma_1,\sigma_2,\mu, \nu$ of degree one  
such that, up to a possible replacement of $A$ by $D$ and of $C$ by $B$, either
\begin{align}  \la{rp11}  &A=\nu \circ z^sR^n(z) \circ \sigma_1^{-1}, & &C=
\sigma_1 \circ z^n \circ \mu \\ 
&\la{rp11+} D=\nu \circ z^n \circ \sigma_2^{-1},& 
&B=\sigma_2 \circ
z^sR(z^n) \circ \mu, \end{align} 
where $R$ is a polynomial, $n\geq 1,$ $s\geq 0$, and $\GCD(s,n)=1,$ or
\begin{align} \la{rp21} &A=\nu \circ T_m \circ \sigma_1^{-1}, &
&C=\sigma_1 \circ T_n \circ \mu, \\
&\la{rp21+} D=\nu \circ T_n \circ \sigma_2^{-1}& &B=\sigma_2 \circ
T_m\circ \mu, \end{align} 
where $T_n, T_m$ are the Chebyshev polynomials, $n,m\geq 1$, and $\GCD(n,m)=1.$ 
\et

Notice that the main difficulty in the practical use of Theorem \ref{rr1} and Theorem \ref{r2}  is the fact that classes of solutions appearing in Theorem \ref{r2} are not disjoint.
Namely, any solution of the form \eqref{rp21}, \eqref{rp21+} with $n=2$
also can be represented in the form \eqref{rp11}, \eqref{rp11+} 
(see e. g. \cite{mz}, \cite{pakmom}, \cite{ms} for further details).

The description of polynomial solutions of equation \eqref{ura} in the general case in a certain sense reduces to the case where \eqref{ggccdd} holds
by the following statement. 

\bt [\cite{en}]\la{r1}
Let $A,C,D,B$ be polynomials 
such that \eqref{ura} holds. Then there exist polynomials
$U, V, \widetilde A, \widetilde C, \widetilde D, \widetilde B,$ where
$$\deg U=\GCD(\deg A,\deg D),  \ \ \ \deg V=\GCD(\deg C,\deg B),$$
such that
$$A=U\circ \widetilde A, \ \  D=U\circ \widetilde D, \ \ C=\widetilde C\circ V, \ \  B=\widetilde B\circ V,$$
and 
$$ \widetilde A\circ \widetilde C=\widetilde D\circ \widetilde B.$$
In particular, if $\deg C=\deg B,$ then 
there exists a polynomial $\mu$ of degree one such that 
$$A=D\circ \mu^{-1}, \ \ \ C=\mu\circ B.$$

\et
\vskip 0.2cm

Theorem \ref{r2} implies  the following 
description of polynomial solutions 
of equation \eqref{1} under the condition \be \la{asdf} \GCD(\deg X, \deg B)=1 \ee  (see \cite{i}). 

\bt [\cite{i}]\la{i}
Let $A,B,X$ be polynomials 
such that \eqref{1} and \eqref{asdf} hold.
Then there exist polynomials $\mu, \nu$ of degree one  
such that either
$$A=\nu \circ z^sR^n(z) \circ \nu^{-1}, \ \  \ X=
\nu \circ z^n \circ \mu, \ \ \ D=\mu^{-1}\circ 
z^sR(z^n) \circ \mu, $$
where $R$ is a polynomial, $n\geq 1,$ $s\geq 0$, and $\GCD(s,n)=1,$ or
$$A=\nu \circ \pm T_m \circ \nu^{-1}, \ \  \ X=
\nu \circ T_n \circ \mu, \ \ \ D=\mu^{-1}\circ 
\pm T_m \circ \mu , $$
where $T_n, T_m$ are the Chebyshev polynomials, $n,m\geq 1$, and $\GCD(n,m)=1.$ 

\et

Notice, however, that Theorem \ref{r2}, even combined with Theorem \ref{r1}, provides very little information 
about solutions of \eqref{1} if \eqref{asdf} is not satisfied. 
A possible way to investigate the general case is to analyze somehow the totality of all decompositions of a polynomial $P$, basing on Theorem \ref{rr1} and Theorem \ref{r2}, and then to apply this analysis to \eqref{1} using the fact that we can pass from the decomposition $P=A\circ X$ to the decomposition $P=X\circ B$.
This way was used in the paper \cite{ms}. A similar techniques was used in the paper \cite{mz} where it was applied to the study 
of decompositions of iterations of a polynomial. In this paper we use another method completely bypassing Theorem \ref{rr1}. 
Notice by the way that Theorem \ref{rr1} does not hold for arbitrary rational functions (see e. g. \cite{mp1}).

\end{section}

\begin{section}{Semiconjugacies and Julia sets}

\begin{subsection}{Polynomials sharing preimages of compact sets}
Let $f_1(z),$ $f_2(z)$ be non-constant complex polynomials
and $K_1,K_2\subset \C$ compact sets. 
In the paper \cite{p1} we investigated  
the following problem.
Under what conditions on the collection $f_1(z), f_2(z), K_1, K_2$ do
the preimages $f_1^{-1}(K_1)$ and $f_2^{-1}(K_2)$
coincide, that is, \be \la{11} f_1^{-1}(K_1)=f_2^{-1}(K_2) =K\ee
for some compact set $K\subset \C$ ?
 
Using ideas from approximation theory,  we relate equation \eqref{11} to the functional equation
\be \l{2}
g_1(f_1(z))=g_2(f_2(z)), \ee
where $f_1(z),f_2(z),g_1(z),g_2(z)$ are polynomials. 
It is easy to see that for any polynomial
solution of \eqref{2} and any compact set $K_3\subset \C$
we obtain a solution of \eqref{11} setting \be \la{3} K_1=g_1^{-1}(K_3), \ \ \
K_2=g_2^{-1}(K_3). \ee Briefly, the main result of  \cite{p1}
states that, under a very mild condition
on the cardinality of $K,$
all solutions of \eqref{11} can be obtained in this way. Combined with Theorem \ref{r1} and Theorem \ref{r2} this leads to a very 
explicit description of  solutions of \eqref{11}.

\bt  [\cite{p1}]\l{zz} Let $f_1(z),$ $f_2(z)$ be
polynomials, $\deg f_1=d_1,$ $\deg f_2=d_2,$ \linebreak $d_1\leq d_2,$ and let $K_1,K_2, K\subset \C$ be compact sets
such that \eqref{11} holds. 
Suppose that
$\rm{card}\{K\}\geq \LCM(d_1,d_2).$ Then, if $d_1$ divides $d_2,$
there exists a polynomial $g_1(z)$ such that $f_2(z)=g_1(f_1(z))$
and $K_1=g_1^{-1}(K_2).$ On the other hand, if $d_1$ does not divide $d_2,$ then there exist polynomials $g_1(z),$ $g_2(z),$
$\deg g_1=d_2/d,$ $\deg g_2=d_1/d,$ where $d=\GCD(d_1,d_2),$ and a compact set $K_3\subset \C$
such that \eqref{2},\eqref{3} hold. Furthermore, in this case there
exist polynomials
$\tilde f_1(z),$ $\tilde f_2(z),$ $W(z),$ $\deg W(z)=d,$
such that 
\be \la{kom} f_1(z)=\tilde f_1(W(z)), \ \ \ \ f_2(z)=\tilde f_2(W(z))\ \ \ \   
\ee 
and 
there exist linear functions $\sigma_1(z),$ $\sigma_2(z)$ such that 
either
\begin{align} \la{t2}  & g_1(z)=z^cR^{d_1/d}(z) \circ \sigma_1^{-1},
&& \tilde f_1(z)  =\sigma_1 \circ z^{d_1/d}, \\
& g_2(z)=z^{d_1/d} \circ \sigma_2^{-1},  && \tilde f_2(z)  =\sigma_2 \circ z^cR(z^{d_1/d}), 
\notag \end{align}
for some 
polynomial $R(z)$ and $c$ 
equal to the remainder after division of $d_2/d$ by $d_1/d,$
or
\begin{align} \la{t3} \ \ \ \ \ \  & g_1(z)= T_{d_2/d}(z)\circ \sigma_1^{-1}, \ \ \ \ \ \ \ \ \ && \tilde f_1(z)  
=\sigma_1 \circ T_{d_1/d}(z),  &&   \\
& g_2(z)= T_{d_1/d}(z)\circ \sigma_2^{-1}, \ \ \ \ \ \ \ \ \ &&  \tilde f_2(z)  =\sigma_2 \circ T_{d_2/d}(z),
\notag\end{align}
for the Chebyshev polynomials $T_{d_1/d}(z),$ $T_{d_2/d}(z).$
\et

Theorem \ref{zz} may be used for proving many other results  
(see \cite{p1} for details) the most notable of which is the following description 
of solutions of \eqref{11} in the case where $K_1=K_2$,
first obtained by T. Dinh (\cite {d}, \cite{d2}) by methods of complex dynamics.

\bt 
[\cite{d2}, \cite{p1}]\l{zzz}
Let $f_1(z),$ $f_2(z)$ be
polynomials such that
\be \la{u22} f_1^{-1}(T)=f_2^{-1}(T)=K \ee
holds for some infinite compact sets  
$T, K\subset \C.$ Then, if $d_1$ divides $d_2,$
there exists a polynomial
$g_1(z)$ such that $f_2(z)=g_1(f_1(z))$ and $g_1^{-1}(T)=T.$
On the other hand, if $d_1$ does not divide $d_2,$ then
there exist polynomials
$\tilde f_1(z),$ $\tilde f_2(z),$ $W(z),$ $\deg W(z)=d,$ 
satisfying \eqref{kom}. 
Furthermore, in this case one of the following conditions holds.
\vskip 0.2cm \noindent
1) $T$ is a union of circles with the common center and 
\be \la{tt2} \tilde f_1(z)=\sigma \circ z^{d_1/d}, \ \ \ \ \ \tilde f_2(z)=\sigma \circ \gamma z^{d_2/d}
\ee for some linear function $\sigma(z)$ and $\gamma \in \C.$ 
\vskip 0.2cm \noindent
2) $T$ is a segment and 
\be \la{tt3} \tilde f_1(z)=\sigma \circ \pm T_{d_1/d}(z), \ \ \ \ \ \tilde f_2(z)=\sigma \circ \pm T_{d_2/d}(z),\ee
for some linear function $\sigma(z)$
and the Chebyshev polynomials $T_{d_1/d}(z),T_{d_2/d}(z).$
\vskip 0.2cm
\et

\end{subsection}

\begin{subsection}{\l{subs} Proofs of Theorem \ref{j+} and Theorem \ref{resh}.}
\noindent{\it Proof of Theorem \ref{j+}.} If $A\underset{X}{\leq} B$, then for any $n\geq 1$ the equality
$$A^{\circ n}\circ X=X\circ B^{\circ n}$$ holds. Therefore,  
if $z_1=X(z_0)$, then the sequence $A^{\circ n}(z_1)$ is bounded if and only if 
the sequence $X\circ B^{\circ n}(z_0)$ is bounded. In  turn, the last sequence is 
bounded if and only if the sequence  $B^{\circ n}(z_0)$ is bounded. Thus,  $A\underset{X}{\leq} B$ implies  \be \l{beg} X^{-1}(K(A))=K(B).\ee 

In other direction, if \eqref{beg} holds, then it follows from $B^{-1}(K(B))=K(B)$ that 
$$(X\circ B)^{-1}(K(A))=K(B).$$ Thus,
$$ X^{-1}(K(A))=(X\circ B)^{-1}(K(A)).$$ 
Since $\deg X\mid \deg (X\circ B)$,  applying to this equality Theorem \ref{zz}   we conclude that 
$$\t A\circ X=X\circ B$$  
for some polynomial $\t A.$ Furthermore, since we proved  that for such $\t A$ the equality  $X^{-1}(K(\t A))=K(B)$ holds, we see that 
$X^{-1}(K(\t A))=X^{-1}(K(A))$, implying that $K(\t A)=K(A)$.
Finally, it follows from Theorem \ref{zz} applied to the equality
$$A^{-1}(K)=\t A^{-1}(K)=K,$$ where $K=K(\t A)=K(A)$,  that 
there exists a polynomial  of degree one $\mu$ such that $\t A=\mu \circ A$ and $\mu(K(A))=K(A).$

More generally, if \be \l{ujn} X^{-1}(K)=K(B)\ee for some compact set $K\subset \C,$
then
$$ X^{-1}(K)=(X\circ B)^{-1}(K),$$ implying by Theorem \ref{zz} that  
equality \eqref{1} holds for some polynomial $A.$ Furthermore, since for such a polynomial  $A$ equality \eqref{beg} holds, we conclude that 
$X^{-1}(K)=X^{-1}(K(A))$ and $K=K(A).$ 
\qed

\bc \l{c3} Let $B$ be a polynomial of degree at least two. Then a polynomial $X$ is contained in $\f E(B)$ if and only $K(B)$ is a union of fibers of $X.$
In particular, if $B_1$ and $B_2$ are polynomials such  that $K(B_1)=K(B_2)$, then $\f E(B_1)=\f E(B_2).$

\ec
\pr Clearly, condition \eqref{ujn} implies that $K(B)$ is a union of fibers of $X.$ In the other direction, if $K(B)$ is a union of fibers of $X,$
then $$K(B)=X^{-1}(X(K(B))),$$ implying that \eqref{ujn} holds for the compact set 
$K=X(K(B))$. \qed

\bc \l{c1}  Let $A,B,$ and $X$ be polynomials such that $A\underset{X}{\leq} B$.
Then for any decomposition $X=X_1\circ X_2$ there exists 
a polynomial $C$ such that 
$$A\underset{X_1}{\leq} C, \ \ \ \ C\underset{X_2}{\leq} B.$$
\ec 
\pr 
By Theorem \ref{j+}, 
$K(B)=X^{-1}(K(A))$. Since $X=X_1\circ X_2$, this implies that  
$K(B)=X_2^{-1}(\t K),$ where $\t K=X_1^{-1}(K(A)).$ Therefore, by Theorem \ref{j+}, there exists  a polynomial $C$ such that 
\be \l{krot} C\circ X_2=X_2\circ B.\ee 
Now 
we have:
$$A\circ X_1\circ X_2= X_1\circ X_2\circ B=X_1\circ C\circ X_2,$$ implying that $A\circ X_1=X_1\circ C$.
\qed


\br Corollary \ref{c1} may be proved  without using Theorem \ref{j+}. Indeed, 
if $X=X_1\circ X_2$, then it follows from the equality 
$$A\circ (X_1\circ X_2)=X_1\circ (X_2\circ B)$$ by Theorem \ref{r1} that 
\be \la{xuli} X_1\circ X_2=U\circ \t W, \ \ \ \ X_2\circ B=V\circ \t W,\ee
where $$\deg \t W=\GCD(\deg (X_1\circ X_2), \deg (X_2\circ B)).$$
Since $\deg X_2\mid \deg \t W$, Theorem \ref{r1} applied to the first equality in \eqref{xuli} 
implies that $\t W=S\circ X_2$ for some polynomial $S$. 
Therefore, 
$$X_2\circ B=V\circ \t W=V\circ S\circ X_2$$
and hence \eqref{krot} holds for $C=V\circ S.$ 

\er

\noindent{\it Proof of Theorem \ref{resh}.} 
By Theorem \ref{j+}, the condition $X_1,X_2\in \f E(B)$ implies that there exist $K_1$, $K_2\subset \C$ such that 
$$X_1^{-1}(K_1)=K(B),\ \ \ \ X_2^{-1}(K_2)=K(B).$$ 
It follows now from Theorem \ref{zz} that  
there exist polynomials 
$X$, $W$, $U_1,$ $U_2$, $V_1,$ $V_2$ such  that $$\deg X=\LCM(\deg X_1,\deg X_2), \ \ \ 
\deg W=\GCD(\deg X_1,\deg X_2),$$
 and equalities 
 $$ X=U_1\circ X_1=U_2\circ X_2$$  and 
\be \l{ddcc} X_1=V_1\circ W, \ \ \ X_2=V_2\circ W\ee
hold. 
Furthermore, there exists $K_3\subset \C$ such that 
$$K_1=U_1^{-1}(K_3), \ \ \
K_2=U_2^{-1}(K_3).$$
Therefore, 
$$X^{-1}(K_3)=K(B),$$ implying by  
Theorem \ref{j+} that $X\in \f E(B)$. Finally, any of equalities \eqref{ddcc} 
implies that  $W\in \f E(B)$ by Corollary \ref{c1}.
 \qed

\end{subsection}

\end{section}

\begin{section}{Semiconjugacies between fixed $A$ and $B$}

\begin{subsection}{\la{appp} \hskip -0.2cm Semiconjugacies between special polynomials} 

For a polynomial $P$ and a finite set $K\subset \C$ denote by $P_{odd}^{-1}(K)$ the subset of $P^{-1}(K)$ consisting of points where the local 
multiplicity of $P$ is odd. Notice that the chain rule implies that if $P=A\circ B$, then 
\be \l{sad} P_{odd}^{-1}(K)=B_{odd}^{-1}(A_{odd}^{-1}(K)).\ee

\bl \la{che} Let $P$ be a polynomial of degree $n\geq 2$, and $K\subset \C$ a finite set containing at least two points. Assume that $P_{odd}^{-1}(K)=K$. Then $K$ contains exactly two points, and $P$ is conjugated to $\pm T_n$.
\el
\pr Denote by $e_z$ the multiplicity of $P$ at $z\in \C,$ and set $r=\mathrm{card}(K)$.
Since for any $y\in \C$ the set $P^{-1}(y)$ contains 
$$n-\sum_{\substack{z\in \C \\ P(z)=y}}(e_z-1)$$ points and 
$$\sum_{z\in \C}(e_z-1)=n-1,$$ we have: \be \la{nerr} 	\mathrm{card}(P^{-1}(K))\geq rn-\sum_{z\in \C}(e_z-1)= (r-1)n+1\ee  (the minimum is attained if $K$ contains all finite critical values of $P$). Therefore, if $$\mathrm{card}(P_{odd}^{-1}(K))=\mathrm{card}(K)=r,$$  then the set $P^{-1}(K)$ contains 
at least $$(r-1)n+1-r$$ points where the local multiplicity of $P$ is greater than one, 
implying that 
\be \la{zxc} \sum_{z\in P^{-1}(K)}\hskip -0.4cm e_z\geq r+2\left((r-1)n+1-r\right).\ee
Since the sum in the left part of \eqref{zxc} equals $rn$, this inequality implies that 
\be \la{pizd} (n-1)(r-2)\leq 0.\ee Thus, $r=2$. Furthermore, since the equality in 
\eqref{pizd} is attained if and only if the equality in \eqref{zxc} is attained, we conclude that if $P_{odd}^{-1}(K)=K$, then 
$e_z=2$ for each  $z\in P^{-1}(K)\setminus K$,  
and the local multiplicity of  $P$ at each of the two points of $K$ is equal to one.

Changing $P$ to $\sigma^{-1}\circ P\circ \sigma$ for a convenient polynomial of degree one $\sigma$,
we can assume that $K=\{-1,1\}.$ 
Then the condition on multiplicities of $P$ implies that $P^2-1$ is divisible by $(P^{\prime})^2$, and calculating the quotient we conclude that $P$ satisfies the differential equation 
$$n^2(1-y^2)=(y^{\prime})^2(1-z^2).$$
Since the general solution of the equation 
$$\frac{y^{\prime}}{\sqrt{1-y^2}}=\pm\frac{n}{\sqrt{1-z^2}}$$ is 
$$\arccos y=\pm n\arccos z+c,$$
it follows now from $P(1)=\pm 1$ that $$P=\pm \cos(n \arccos x)=\pm T_n(z).\eqno\qed$$ 

\br \l{sopr} Notice that the equality $T_n(-z)=(-1)^nT_n(z)$ implies that for even $n$ the polynomials 
$T_n$ and  $-T_n$ are conjugated  since  $T_n=\alpha \circ (-T_n)\circ \alpha^{-1},$ where $\alpha(z)=-z.$
For odd $n$ however the polynomials $T_n$ and  $-T_n$ are not conjugated.
\er

\bl \la{proo} Let $P$ be a polynomial and $a,b\in \C.$ Then the set $P_{odd}^{-1}\{a,b\}$ contains 
at least two points.
\el
\pr It follows from the equality $$2n=\sum_{\substack{z\in \C \\ P(z)=a}}e_z
+\sum_{\substack{z\in \C \\ P(z)=b}}e_z$$ that the number
$$\sum_{z\in P_{odd}^{-1}\{a,b\}}e_z$$ is  even, implying that the number
$\mathrm{card}(P_{odd}^{-1}\{a,b\})$ also is even. On the other hand, 
$$\mathrm{card}(P_{odd}^{-1}\{a,b\})\neq 0,$$ 
for otherwise  $P_{odd}^{-1}\{a,b\}$ contains at most $n/2+n/2=n$ points in contradiction 
with inequality \eqref{nerr}.
\qed


\bt \la{sem} Let $A$ and $B$ be polynomials of degree at least two  such that \linebreak $A\leq B.$ Then 
$A$ is conjugated to $z^n$ if and only if $B$ is conjugated to $z^n$. 
Similarly, $A$ is conjugated to $\pm T_n$ if and only if $B$ is conjugated to $\pm T_n.$
\et

\pr 
Assume that $B$ is conjugated to $\pm T_n$, and let $X$ be a semiconjugacy from $B$ to $A$. 
Changing $B$ and $X$ to $\sigma^{-1}\circ B\circ \sigma$ and 
$X\circ\sigma$, 
for a convenient polynomial  $\sigma$ of degree one, without loss of generality we can assume that 
$B=\pm T_n.$ By Theorem \ref{j+}, we have:
\be \l{korova} X^{-1}(K(A))=K(B)=[-1,1].\ee
Set $m=\deg X.$ Since  \be \l{smer} T_{m}^{-1}([-1,1])=[-1,1],\ee equality
\eqref{korova} implies that   
$$X^{-1}(K(A))=T_{m}^{-1}([-1,1]).$$
It follows now from Theorem \ref{zz} that there exists a polynomial  $\delta$ of degree one such that 
$X=\delta \circ T_m.$ Therefore, changing  $A$ and $X$ to $\delta^{-1}\circ A\circ \delta$ and 
$\sigma^{-1}\circ X$, we can assume that $X=T_m$.  
Thus, we have: \be \l{pesa} A\circ  T_m= T_{m}\circ \pm T_n=(-1)^m T_n\circ T_m,\ee implying that
$A=\pm T_n$.

Similarly, if $B=z^n$, then the equalities 
$$ X^{-1}(K(A))=K(B)=\mathbb D,$$
and $(z^m)^{-1}(\mathbb D)=\mathbb D$ imply that $X=\delta \circ z^m$ for some polynomial  $\delta$ of degree one, 
and arguing as above we conclude that $A$ 
is conjugated to $z^n.$

Assume now that $A$ is conjugated to $\pm T_n$. 
Without loss of generality we can assume that $A=\pm T_n.$
Since 
$T_{n\, odd}^{-1}\{-1,1\}=\{-1,1\}$, formula \eqref{sad} implies that 
 $$(\pm T_n\circ X)_{odd}^{-1}\{-1,1\} =X_{odd}^{-1}\{-1,1\}.$$  It follows now from  
\be \l{qazz} \pm T_n\circ X=X\circ B\ee that
\be \la{bb} B_{odd}^{-1}(X_{odd}^{-1}\{-1,1\})=X_{odd}^{-1}\{-1,1\}.\ee  
Since by Lemma \ref{proo} the set $X_{odd}^{-1}\{-1,1\}$ contains at least two points, this implies by   
Lemma \ref{che} that the polynomial $B$ is conjugated to $\pm T_n.$  

Finally, if $A$ is conjugated to $z^n$, 
we can assume that $A=z^n$, and considering zeroes of the left and the right parts of the equality
$$z^n\circ X=X\circ B,$$ we see
that $B^{-1}(X^{-1}(0))=X^{-1}(0).$ It follows now from inequality \eqref{nerr} that $X^{-1}(0)$ consists of a single point, implying easily that the polynomial $B$ is conjugated to $z^n.$ 
\qed

\br \l{sopr2} Since for even $n$  the polynomials $T_n$ and $-T_n$ are conjugated (see Remark \ref{sopr}), Theorem \ref{sem} implies that if
$B$ is conjugated to $\pm T_n$ for even $n$, then $A$ and   $B$ are  conjugated. On the other hand, 
if $B$ is conjugated to $-T_n$ for odd $n$,  then $A$ is not necessary conjugated to $-T_n$, but only to $\pm T_n$.
Still, it follows from \eqref{pesa} that  
if  $B$ is conjugated to $T_n$,
then $A$ is conjugated to $T_n$. 
\er 

Notice that Theorem \ref{sem} combined  with Remark \ref{sopr2} implies the following corollary. 

\bc  \l{kote}
Let  $A$ and $B$ be polynomials such that the conditions $A\leq B$ and $B\leq A$ hold simultaneously, and at least one of $A$ and $B$ is special.
Then $A$ and $B$ are conjugated. \qed
\ec

\end{subsection}

\begin{subsection}{Proof of Theorem \ref{uni}}
The following lemma is a well-known fact from the complex dynamics. For the reader's convenience we 
give a short proof based on Theorem \ref{zz}.

\bl \l{jp} Let $A$ be a polynomial of degree $n$ such that $K(A)$ is a union of circles with a common center. Then $K(A)$ is a disk, and $A$ is conjugate to $z^n.$
Similarly, if  $K(A)$ is a segment, then $A$ is conjugated to $\pm T_n.$ 
\el
\pr Since for a polynomial $A$ the complement to $K(A)$ in $\C\P^1$ is connected (see e.g. \cite{mil}, Lemma 9.4), if $K(A)$ is a union of circles with a common center, then  $K(A)$ is a disk. Furthermore, 
changing if necessary $A$ to a conjugated polynomial, we can assume that $K(A)=\mathbb D.$ Thus,  
$A^{-1}(\mathbb D)=\mathbb D.$
On the other hand, $(z^n)^{-1}(\mathbb D)=\mathbb D,$ and applying to these equalities   Theorem \ref{zz}, we conclude that
$ A=\alpha z^n,$ where  $\vert \alpha\vert =1,$ implying that
$A$ is conjugate to $z^n.$

Similarly, if $K(A)$ is a segment, we can assume that $K(A)=[-1,1]$, and to conclude in a similar  way 
that $A$ is conjugated to $\pm T_n.$
\qed
\vskip 0.2cm
\noindent{\it Proof of Theorem \ref{uni}.}
Set $d_0=\deg X_0$, and let $X\in\f E(A,B)$ be a polynomial of degree $d$. 
By Theorem \ref{j+}, we have: 
$$ X_0^{-1}(K(A))=K(B), \ \ \ \   X^{-1}(K(A))=K(B).$$
Applying to these equalities Theorem \ref{zzz}
and taking into account that, by Lemma \ref{jp},  $K(A)$ is neither a union of circles with the common center nor a segment,
we conclude that  
$X=\t A\circ X_0$ for some polynomial $\t A$. Substituting now this expression in \eqref{1} and using that $X_0\in\f E(A,B)$ we have:  
$$A\circ \t A\circ X_0=\t A\circ X_0\circ B=\t A\circ A\circ X_0,$$ implying that 
$A\circ \t A=A\circ \t A$. 

In the other direction, if $A$ commutes with $\t A$, then 
$$ A\circ (\t A\circ X_0)=\t A\circ A\circ X_0=(\t A\circ X_0)\circ B. \eqno{\Box}$$

\vskip 0.2cm
 Theorem \ref{uni} implies in particular the following classification of commuting polynomials obtained by Ritt.

\bt [\cite{r}]\la{commut} Let $A$ and $B$ be commuting polynomials of degree at least two. Then, up to the change
\be \l{c} A\rightarrow \lambda \circ A\circ  \lambda^{-1},\ \ \ B\rightarrow \lambda \circ B\circ  \lambda^{-1},\ee where $\lambda$ is a polynomial of degree one, either 
\be \l{cc1} A=z^n, \ \ \ B=\v z^m,\ee where $\v^n=\v,$ or
\be \l{cc2} A=\pm T_n, \ \ \ B=\pm T_m,\ee
or 
\be \l{iri} A= \v_1R^{\circ m}, \ \ \ \ B=\v_2R^{\circ n},\ee
where $R=zS(z^{\ell})$ for some polynomial $S$, and  $\v_1,$ $\v_2$ are $l$-th roots of unity. 
\et
\pr Assume first that $A$ is conjugated to $z^n$. Without loss of generality we may assume that $A=z^n.$ 
Applying Theorem \ref{j+} for $B=A$ and $X=B$,   we have:
$$B^{-1}(K(A))=K(A).$$ Since $K(A)=\mathbb D$, arguing as in Lemma \ref{jp} we conclude that  
$B=\v z^m$, and it follows from $A\circ B=B\circ A$ that  $\v^n =\v.$
If $A$ is conjugated to $\pm T_n$, the proof is similar.

On the other hand, if $A$ is non-special, then  Theorem \ref{uni} implies that any $B\in \f E(A,A)$ has the form 
$B=\t A\circ R$, where $R$ is a polynomial of the minimum possible degree in $\f E(A,A)$.
Now we can apply Theorem \ref{uni} again to the polynomial $\t A$ and so on, arriving eventually to the representation 
$B=\mu_1 \circ R^{\circ m_1}$, where $\mu_1$ is a polynomial  of degree one commuting with $A$. 
In particular, since $A\in \f E(A,A)$, the equality  
$A=\mu_2\circ R^{\circ m_2}$ holds 
for some polynomial  $\mu_2$ of degree one  commuting with $A$.
Furthermore, since $R$ commutes with $A=\mu_2\circ R^{\circ m_2}$, the
polynomial $\mu_2$ commutes with $R$. This implies easily that, up to a conjugacy, 
$R=zS(z^{\ell})$ for some polynomial $S$, and $\mu_2=\v_2z$ for some $l$th root of unity $\v_2.$
Finally, since $\mu_1$ commutes with the polynomial $A$, and $A=\mu_2\circ R^{\circ m_2}$ has the form $z\t S(z^{\ell})$ for some polynomial $\t S$, 
we conclude that $\mu_1=\v_1 z$ for some $l$th root of unity $\v_1.$  \qed

\end{subsection}

\begin{subsection}{\la{mesc} Semiconjugacies and invariant curves}
It was shown in the recent paper \cite{ms} that the problem of describing of 
semiconjugate polynomials is closely related to the  problem of describing of algebraic
curves $\f C$ in $\C^2$ invariant under maps of the form $F:\,(x,y)\rightarrow (f(x),g(y)),$ where $f,g$ 
are polynomials of degree at least two.
Briefly, this relation may be summarized as follows (see Proposition 2.34
of \cite{ms} for more details).

If $\f C$ is an irreducible $(f,g)$-invariant curve, then its projective closure $\overline {\f C}$ in $\C\P ^1\times \C\P^1$ is also $(f,g)$-invariant. Denote by
$\bar h$ the restriction  of $F$ on  $\overline {\f C}$.
Let $\t {\f C}$ be the desingularization of $\f C$ and $\beta\,: \widetilde {\f C}\rightarrow \overline {\f C}$ a map biholomorphic off a finite set.
Clearly, $\bar h$
lifts to a holomorphic map $h\,:\widetilde {\f C}\rightarrow \widetilde {\f C}.$ 
Consider now the commutative diagram 
\be 
\begin{CD} \l{opa}
\widetilde {\f C} @>h>> \widetilde {\f C}\\
@VV \beta V @VV \beta V\\ 
 \overline {\f C} @>\bar h>>\overline {\f C} \\
@VV \alpha V @VV \alpha V\\ 
\C\P^1 @>f >> \C\P^1 ,   
\end{CD}
\ee
where 
$\alpha\,:\overline {\f C}\rightarrow \C\P^1$ is the projection map  
onto the first coordinate.  Set $\pi=\alpha\circ \beta$. If $\pi$ is a constant, 
then $\f C$ is a line $z_1=\xi,$ where $\xi$ is a fixed point of $f$, so assume that 
the degree of $\pi$ is at least one.
Observe that since $f^{-1}(\infty)=\infty$, 
the set $K= \pi^{-1}(\infty)$ and the map $h$ satisfy the equality   
\be \l{xomiak} h^{-1}(K)=K.\ee

Since $h$ is a holomorphic map between Riemann surfaces of the same genus and $\deg h=\deg f\geq 2$,
it follows from the Riemann-Hurwitz formula that either $g(\t {\f C})=0$, or $g(\t {\f C})=1$ and 
$h$ is unbranched. Since $\deg h\geq 2 $,  for unbranched $h$ equality \eqref{xomiak} is impossible. Therefore, 
$\widetilde {\f C}=\C\P^1$ and \eqref{xomiak} implies easily that, up to the change $\alpha\circ h \circ \alpha^{-1},$ where $\alpha$ is a M\"obius transformation, either $K=\infty$ and 
$h$ is a polynomial, or $K=\{0,\infty\}$ and $h=z^{\pm\deg f}.$ Thus, 
\be \la{pl1} f\circ \pi =\pi \circ h,\ee where either $\pi$ and $h$ are polynomials, or 
$h=z^{\pm \deg f}$ and $\pi $ is a Laurent polynomial. The last case requires an additional investigation. The paper \cite{ms}  refers 
(Fact 2.25) to a more general  result of \cite{medv} (Theorem 10) implying that for a non-special polynomial $f$
this possibility is excluded. Alternatively, one can use results of the paper \cite{pak} 
(e.g. Theorem 6.4).

Considering in a similar way the projection onto 
the second coordinate, we arrive to the equality 
\be \la{pl2} g\circ \rho =\rho \circ h.\ee Thus, for non-special $f$ and $g$ any irreducible 
$(f,g)$-invariant curve may be paramet\-rized by some polynomials $\pi,$ $\rho$ satisfying a 
system given by equations \eqref{pl1}, \eqref{pl2} for some polynomial $h$.

Notice that in a certain sense a description of $(f,g)$-invariant curves 
reduces to the case $f=g$ since the commutative diagram 
\be
\begin{CD} 
 \C^2 @>(h,h)>> \C^2\\
@VV (\pi,\rho) V @VV (\pi,\rho) V\\ 
\C^2 @>(f,g)>>\C^2 
\end{CD}
\ee
implies that any $(f,g)$-invariant curve is an image of an $(h,h)$-invariant curve 
under the map $(x,y)\rightarrow (\pi(x),\rho(y)).$

Theorem \ref{uni} permits to obtain easily the following description of $(f,f)$-invariant curves obtained in \cite{ms} (see Theorem 6.24 and the theorem on p. 85).

\bt Let $f$ be a non-special polynomial of degree at least two, and  $\f C$ an irreducible $(f,f)$-invariant curve in $\C^2$. Then there exists  a polynomial $p$ which commutes with $f$ such that $\f C$ has either 
 the form $z_1= p(z_2)$ or $z_2= p(z_1)$.
\et
\pr If $\f C$ is a line $z_1=\xi,$ then $\xi$ is a fixed point of $f$,
and the conclusion of the theorem holds for $p=\xi$. Similarly, the theorem holds if $\f C$ is a line $z_2=\xi.$ Otherwise, as it was shown above, $\f C$ may be paramet\-rized by some non-constant polynomials $\pi,$ $\rho$ satisfying the system 
\be \la{pl1+} f\circ \pi =\pi \circ h,\ee
\be \la{pl2+} f\circ \rho =\rho \circ h\ee
 for some polynomial $h$.
Furthermore, without loss of generality we may assume that there exists no polynomial $w$ of degree greater than one such that 
\be \l{emat} \pi =\t \pi\circ w, \ \ \ \rho =\t \rho\circ w\ee 
for some polynomials $\t \pi,$ $\t \rho$. Indeed, if \eqref{emat} holds, then applying Theorem \ref{r1} to the 
equality $$(f\circ \t \pi)\circ w= \t \pi\circ (w\circ h),$$  
we conclude that $w\circ h=\t h \circ w$ for some polynomial $\t h$, implying that 
we may change $\pi$ to $\t \pi,$ $\rho$ to $\t \rho,$ and $h$ to $\t h$. 


Set $d=\GCD(\deg \rho, \deg \pi)$. Since $f$ is not special, it follows from \eqref{pl1+}, \eqref{pl2+} by Theorem \ref{uni} that if both $\rho$ and $\pi$ are of degree at least two, then $d>1$, implying by Theorem \ref{resh} that \eqref{emat} holds for some 
 polynomials $\t \pi,$ $\t \rho$ and $w$ with $\deg w=d>1.$ Therefore, at least one of 
polynomial $\rho$ and $\tau$ is of degree one. Assume say that   
$\deg \rho= 1$. Then, $\f C$ has the form 
$z_1=p(z_2),$ where $p=\pi\circ \rho^{-1}$. Furthermore, equality \eqref{pl2+} implies that 
$h=\rho^{-1}\circ f \circ \rho,$ and substituting this expression into \eqref{pl1+} we conclude that $p$   commutes with $f.$ \qed

\vskip 0.2cm

\noindent{\it Proof of Theorem \ref{new}.}
For any polynomials of coprime degrees $u$ and $v$ the curve $\f C_{u,v}:\ u(x)-v(y)=0$ is irreducible (see e.g. \cite{pakcur}, Proposition 3.1). Furthermore, if 
\eqref{krys} holds and  $(x_0,y_0)$ is a point on $\f C_{u,v}$, then \eqref{krys} yields the equality  
$$u(f(x_0))=t(u(x_0))=t(v(y_0))=v(g(y_0)),$$ implying that
$(f(x_0),g(y_0))$ also is a point on $\f C_{u,v}$.

In the other direction, assume that $\f C$ is an irreducible $(f,g)$-invariant curve which is not a line, and let $\pi$ and $\rho$ be polynomials  parametrizing $\f C$ and satisfying \eqref{pl1}, \eqref{pl2} for some polynomial $h.$ 
Then by Theorem  \ref{resh}, there exist polynomials $u$ and $v$ of coprime degrees such that 
$$u\circ \pi=v\circ \rho.$$ Thus, any irreducible $(f,g)$-invariant curve $\f C$ in $\C^2$ has
the form $u(x)-v(y)=0$ for some polynomials $u,v$ of coprime degrees.
Furthermore, since the polynomial
$$t=u\circ \pi=v\circ \rho$$ is contained in $\f E(h)$ we have:
$$t\circ u\circ \pi=u\circ \pi\circ h=u\circ f\circ \pi,$$ 
$$t\circ v\circ \rho=v\circ \rho\circ h=v\circ g\circ \rho,$$ implying \eqref{krys}. \qed

\vskip 0.2cm
A further analysis of system \eqref{krys} using Proposition \ref{predl} and Proposition \ref{suks} proved below leads to a more precise description of  $(f,g)$-invariant curves apparently equivalent to the one given by Theorem 6.2 of \cite{ms}. 
Notice however that in the paper \cite{ms} a more general  case of {\it skew-invariant} curves and  {\it skew-twists} between polynomials is considered, and 
the methods of our paper involving Julia sets seem not to be extendable to this more general situation.

\end{subsection}

\begin{subsection}{Semiconjugacies between equivalent $A$ and $B$}

For a natural number $n>1$ with a prime decomposition $n=p_1^{a_1}p_2^{a_2}\dots p_s^{a_s}$
set  $\mathrm{rad}(n)=p_1p_2\dots p_s$. The following two theorems in totality provide  a proof of Theorem \ref{e}.

\bt \l{wsx} 
Let $A$ and $B$ be polynomials of degree at least two. Then conditions $A\leq B$ and $B\leq A$  hold simultaneously
if and only if  $A\sim B$.
\et 
\pr The ``if'' part follows from the definition of $\sim$ (see the introduction). Furthermore, if at least one of $A$ and $B$ is special, then conditions $A\leq B$ and $B\leq A$ 
imply by Corollary \ref{kote} that  $A$ and $B$ are conjugated and hence equivalent. So, we may assume that $A$ and $B$ are non-special.

Let $Y$ and $X$ be 
polynomials such that 
\be \l{vbn} B\underset{Y}{\leq} A, \ \ \ \  A\underset{X}{\leq} B.\ee Set $n=\deg A=\deg B.$ We can assume that $\deg X>1,$ $\deg Y>1$ since otherwise 
$A$ and $B$ are conjugated and hence $A\sim B.$
Since \eqref{vbn} implies that $Y\circ X$ commutes with $B$, Theorem \ref{commut} implies that 
\be \l{divi} \rad(\deg X)\mid \rad(n).\ee
In particular,
\be \la{fgd} \GCD(\deg X,n)>1.\ee 
Applying Theorem \ref{r1} to the equality \be \la{yuy} A\circ X=X\circ B,\ee   
we conclude that there exist polynomials $\tt X,$ $\tt B$, and $W$  such that  
\be \l{of} B=\tt B\circ W, \ \ \ X=\tt X\circ W,\ee 
and $\deg W=\GCD(\deg X,n).$
Clearly, $B\sim  W\circ \t B$, and 
equalities \eqref{yuy} and \eqref{of} imply that \be \la{guk} A\circ \t X=\t X\circ (W\circ \t B).\ee  
Furthermore,
$\deg \t X<\deg X$, since $\deg W>1$ by \eqref{fgd}. 
If $\deg \t X=1,$ then $A\sim W\circ \t B$ since $A$ and $W\circ \t B$ are conjugated;
hence, $$A\sim W\circ \t B\sim B,$$ and we are done. Otherwise, 
we can apply Theorem \ref{r1} in a similar way to equality \eqref{guk} and so on. Since 
condition \eqref{divi} ensures that the degrees of corresponding semiconjugacies decrease,  
we obtain in this way  a finite chain of equivalences from $B$ to $A$. \qed

\bt \l{xriak} Let $A$ and $B$ be polynomials of degree at least two.  Then $A\sim B$
if and only if there exist polynomials $X$ and $Y$  
such that 
\be \la{ioi+} B\circ Y=Y\circ A, \ \ \ \ A\circ X=X\circ B,\ee  
and $Y\circ X=B^{\circ d}$ for some $d \geq 0.$ 
\et 
\pr 
Taking into account Theorem \ref{wsx}, we only must show that if equalities \eqref{ioi+} hold, then they hold for some $\t X,$ $\t Y$ such that 
 $\t Y\circ \t X=B^{\circ d}$, $d \geq 0.$ 
Since \eqref{ioi+} implies that $Y\circ X$ commutes with $B$, it follows from Theorem \ref{commut} that either $B$ is special, or, up to a conjugacy, 
$$Y\circ X= \v_1 R^{\circ m_1}, \ \ \ B=\v_2 R^{\circ m_2},$$ where $R=zS(z^n)$ for some polynomial $S$, and $\v_1,$ $\v_2$ are $n$th roots of unity.
In the first case, Corollary \ref{kote} implies that $A$ and $B$ are conjugated. Therefore, in this case \eqref{ioi+} holds for some M\"obius transformations $\t Y$ and  
$\t X$ such that $\t Y\circ \t X=B^0.$ 
In the second case set $$\t X=X\circ \v_3R^{\circ (m_2m_1-m_1)},$$ where $\v_3=\v_2^{m_1}/\v_1,$ and observe that the second of equalities \eqref{ioi+} still holds for $\t X$ since 
$$A\circ \t X=A\circ X\circ  \v_3R^{\circ (m_2m_1-m_1)}=X\circ B\circ  \v_3R^{\circ (m_2m_1-m_1)}=$$
$$=X\circ \v_2 R^{\circ m_2}\circ  \v_3R^{\circ (m_2m_1-m_1)}=X\circ\v_3R^{\circ (m_2m_1-m_1)}\circ  \v_2 R^{\circ m_2}=\t X\circ B. $$
On the other hand, we have:
$$Y\circ \t X=\v_1 R^{\circ m_1}\circ \v_3R^{\circ (m_2m_1-m_1)}=\v_1 \v_3R^{\circ m_2m_1}= \v_2^{m_1}R^{\circ m_2m_1}=B^{\circ m_1}.\eqno{\Box}$$

\end{subsection}

\end{section}

\begin{section}{Semiconjugacies for fixed  $B$}
\begin{subsection}{Special factors of semiconjugacies}

\bl \la{xyx} Let $A$ and $B$ be polynomials of degree $n\geq 2$ such that 
\be \la{chi} A\circ T_{\ell}=T_{\ell}\circ B, \ \ \ l\geq 2.\ee
Then $l\leq 2n,$ unless $A=\pm T_n$ and  $B=\pm T_n.$ 
Similarly, if
\be \la{ci} A\circ z^{\ell}=z^{\ell}\circ B, \ \ \ l\geq 2,\ee then $l\leq n$, unless $A=\alpha z^n,$ $\alpha \in \C,$
and $B=\beta z^n,$ $\beta \in \C.$
\el
\pr 
If 
\be \la{nera} n\leq \frac{l-1}{2},\ee then the set $$(T_{\ell}\circ B)_{odd}^{-1}\{-1,1\}=B_{odd}^{-1}\{-1,1\}$$ contains at most $l- 1$ points. Therefore, if equality \eqref{chi} holds, then the set 
\be \la{set} (A\circ T_{\ell})_{odd}^{-1}\{-1,1\}\ee also contains at most $l- 1$ points. On the other hand, 
since  $-1$ and $1$ are the only finite critical values of $T_n$, if the set $A_{odd}^{-1}\{-1,1\}$ contains at least one point distinct from $\pm 1$, then set \eqref{set}  
contains at least $l$ points.
Since by Lemma \ref{proo} the  set $A_{odd}^{-1}\{-1,1\}$ contains at least two points, we conclude that if \eqref{nera} holds, then 
\be \l{op} A_{odd}^{-1}\{-1,1\}=\{-1,1\}.\ee Therefore, by Lemma \ref{che}, $A=\pm T_n$,  
It follows now from \eqref{chi} that $$\pm T_{nl}=T_{\ell}\circ B,$$ implying that $$T_{\ell}\circ B=\pm T_{\ell}\circ T_n,$$ and applying to the last equality Theorem \ref{r1} 
we see that 
\be \la{iuy} T_{\ell}=\pm T_{\ell}\circ \mu, \ \  \ B=\mu^{-1}\circ T_{n},\ee for some polynomial   $\mu$ of degree one.  
Finally, it is easy to see, using for example the explicit formula 
\be \la{cheb} T_n=\frac{n}{2}\sum_{k=0}^{[n/2]}(-1)^k\frac{(n-k-1)!}{k!(n-2k)!}(2x)^{n-2k},\ee 
that $T_n$ has non-zero coefficients of its terms of degree $n$ and $n-2$, and the coefficient equal zero for its term of degree $n-1.$ Thus,  
the first of equalities \eqref{iuy} 
implies the equality $\mu =\pm x$.

Assume now that equality \eqref{ci} holds and $n\leq l-1$. Then the polynomial in the right part of \eqref{ci} has at most $l-1$ zeroes.
On the other hand, since the unique finite critical value of $z^{\ell}$ is zero, it is easy to see that, unless \be \la{wer} A=\alpha z^n, \ \ \ \alpha\in \C,\ee
the polynomial in the left part of \eqref{ci} has at least $l$ zeroes. 
Finally, \eqref{wer} and \eqref{ci} imply easily that 
$B=\beta z^n,$ $\beta\in \C.$   
\qed

\bt \la{ui} Let $B$ be a non-special polynomial of degree $n\geq 2$, 
and $X$ an element of $\f E(B).$
Assume that $X=W_1\circ z^{\ell}\circ W_2$ for some polynomials $W_1,$ $W_2$ and $l\geq 1$.
Then $l\leq  n.$  Similarly, if $X=W_1\circ \pm T_{\ell}\circ W_2$,
then $l\leq  2n.$ 
\et 
\pr If $X=W_1\circ z^{\ell}\circ W_2$, then applying Corollary \ref{c1} twice we conclude that there exist polynomials $C_1$and $C_2$ such that the equalities 
\be \l{pok} A\circ W_1=W_1\circ C_1, \ \ \ C_1\circ z^{\ell}=z^{\ell}\circ C_2, \ \ \ C_2\circ W_2=W_2\circ B
\ee
hold. Applying now Lemma \ref{xyx} to the second equality in \eqref{pok} we conclude that $l\leq n$, unless $C_1$ and $C_2$ are conjugated to $z^n.$ On the other hand, in the last case  the third equality in \eqref{pok} implies  by Theorem \ref{sem} that $B$ is conjugated to $z^n.$  
If $X=W_1\circ \pm T_{\ell}\circ W_2$, the proof is similar. \qed

\bc \la{iter} Let $B$ be a non-special polynomial of degree $n\geq 2$.
Assume that $B^{\circ d}=W_1\circ z^{\ell}\circ W_2$ for some polynomials $W_1,$ $W_2$, and $l\geq 1,$ $d\geq 1.$
Then $l\leq  n.$  Similarly, if $B^{\circ d}=W_1\circ \pm T_{\ell}\circ W_2$,
then $l\leq  2n.$ 

\ec
\pr Follows from Theorem \ref{ui}, since $B^{\circ d}$ is a semiconjugacy from $B$ to $B$. \qed

\end{subsection}

\begin{subsection}{Proof of Theorem \ref{exi}}
For natural numbers $n$ and $m$ define $l=l(n,m)$ as the maximum number  coprime with $n$ which divides $m$. Thus,
\be \l{ee} m=lb,\ee where $\mathrm{rad}(b)\vert \mathrm{rad}(n)$ and $\GCD(n,l)=1$. Define now $d=d(n,m)$ as the minimum number such that $b$ in \eqref{ee} satisfies $b\mid n^d$. The next proposition describes a general structure of elements of $\f E(B)$
for non-special $B$.

\bp \l{predl}  Let $B$ be a non-special polynomial of degree $n\geq 2$. Then any $X\in\f E(B)$ has the form  
$X=\nu\circ z^{l(n,m)} \circ W,$ where $\nu$ is a polynomial of degree one, and $W$ is a compositional right factor of $B^{\circ d(n,m)}$. 
Furthermore, 
$l(n,m)<n$. 
\ep 
\pr Set $m=\deg X$, and let $l,b,d$ be the numbers defined above. If $A$ is a polynomial such that \be \l{qqaa} A\circ X=X\circ B,\ee then the equality 
\be \la{skuns} A^{\circ d}\circ X=X\circ B^{\circ d},\ee 
implies by Theorem \ref{r1} that 
\be \la{pred} X=U\circ S, \ \ \ B^{\circ d}=V\circ S,\ee for some polynomials $U,V,S$, where $\deg U=l.$
Furthermore, equalities \eqref{qqaa} and  $X=U\circ S$ imply by Corollary \ref{c1} that 
\be \l{bred} A\circ U=U\circ C\ee for some polynomial $C$. 
Since $l$ is coprime with $n$, 
by Theorem \ref{i}  
there exist polynomials $\mu, \nu$ of degree one  
such that either
$$A=\nu \circ z^sR^{\ell}(z) \circ \nu^{-1}, \ \  \ U=
\nu \circ z^{\ell} \circ \mu, \ \ \ C=\mu^{-1}\circ 
z^sR(z^{\ell}) \circ \mu, $$
where $R$ is a polynomial, $n\geq 1,$ $s\geq 0$, and $\GCD(s,l)=1,$ or
$$A=\nu \circ \pm T_n \circ \nu^{-1}, \ \  \ U=
\nu \circ T_{\ell} \circ \mu, \ \ \ C=\mu^{-1}\circ \pm
T_n \circ \mu, $$
where $\GCD(l,n)=1.$
In the last case however Theorem \ref{sem} applied to \eqref{qqaa} implies that $B$ is conjugated to $T_n.$ 
Therefore, the first case must hold and hence $X=\nu\circ z^{\ell}\circ W$, where $W=\mu\circ S$ is a compositional right factor of $B^{\circ d}.$
Moreover, since $n=rl+s,$ where $r=\deg R,$ the inequality $l<n$ holds whenever $r\neq 0.$ On the other hand, if $r=0$, then  $A$ is conjugated to
$z^n$ and hence $B$ also is conjugated to 
$z^n$ by Theorem \ref{sem}.
\qed
\vskip  0.2cm
For a natural number $n>1$ with a prime decomposition $n=p_1^{a_1}p_2^{a_2}\dots p_s^{a_s}$
set  $\ord_p(n)=a_i$, if $p=p_i$ for some $i$, $1\leq i \leq s,$ and $\ord_pn=0$ otherwise.

\bp \l{suks} If, under assumptions of Proposition \ref{predl}, the polynomial $X$ is not a polynomial in $B$, then 
$d(n,m)\leq 2\log_2n+3.$
\ep
\pr Set \be \l{net} a=n^{d}/b.\ee Clearly, for any prime $p,$
$$ \ord_p(b)+\ord_p(a)= \ord_p(n)d,$$
implying that
\be \l{barsuk} \ord_pb= \ord_p(n)(d-1)+\ord_p(n)-\ord_p(a).\ee 
Observe that the definition of $d(n,m)$ implies that $a$ is not divisible by $n$. Moreover, the number $b$ is not divisible by $n$ either, since otherwise equality \eqref{qqaa} implies 
by Theorem \ref{r1} that $X$ is a polynomial in $B$. 
Observe also  that by Theorem \ref{sem} any  polynomial 
$A$ such that \eqref{qqaa} holds is not special.

It follows from Theorem \ref{r1} applied to  equality \eqref{skuns} that 
there exist polynomials $N$, $F$ and $Y$, $Z$, where $$ \deg Z=l, \ \ \ \ \deg Y=a,$$ such that 
$$A^{\circ d}=N\circ Y, \ \ \ \ X=N\circ Z,$$ and 
\be \l{xe} Y\circ X=Z\circ B^{\circ d}.\ee
Applying now Theorem \ref{r1} and Theorem \ref{r2} to the 
equality 
$$Y\circ X=(Z\circ B^{d-i})\circ B^i$$ for each $i,$ $1\leq i \leq d-1$, we obtain a collection of polynomials  
$Y_i,$ $X_i$, $W_i$ $U_i,$ $K_i$, $L_i,$ $1\leq i \leq d-1,$  
such that 
\be \la{ewq} Y=U_i\circ Y_i, \ \   Z\circ B^{\circ d-i}=U_i\circ K_i, \ \ 
X=X_i\circ W_i, \ \ B^{\circ i}=L_i\circ W_i,\ee and 
\be \l{ebs} Y_i\circ X_i=K_i \circ  L_i.\ee
Furthermore,
$$\deg Y_i=a_i,\ \ \ \deg X_i=lb_i, $$ 
where 
\be \l{brose} a_i=\frac{a}{\GCD(a,n^{d-i})},\ \ \ b_i=\frac{b}{\GCD(b,n^i)},\ee
and 
there exist  polynomials of degree one 
$\nu_i,$ $\sigma_i,$ $\mu_i$  
$1\leq i \leq d-1,$  such that either 

\be \l{en1}Y_i = \nu_i\circ z^{a_i}\circ \sigma_i, \ \ \ \  X_i=\sigma_i^{-1}\circ z^cR(z^{a_{i}})\circ \mu_i, \ee where $R\in \C[z]$ and $\GCD(c,a_i)=1,$ or 

\be \l{en2} Y_i= \nu_i\circ z^cR^{lb_i}(z) \circ \sigma_i, \ \ \ \  X_i =\sigma_i^{-1}\circ z^{lb_i}\circ \mu_i,\ee where $R\in \C[z]$ and $\GCD(c,lb_i)=1,$ or
\be \l{en3} Y_i= \nu_i\circ T_{a_i} \circ \sigma_i, \ \ \ \  X_i =\sigma_i^{-1}\circ T_{lb_i}\circ \mu_i,\ee 
where $\GCD(a_i,lb_i)=1.$ 

Observe first that \be \la{neru} a_i\geq 2^i,\ \ \ b_i\geq 2^{d-i}.\ee 
Indeed, since $n\nmid a$,  there 
exists $p\in \mathrm{rad}(n)$ such that $\ord_p(n)-\ord_p(a)>0$. Thus, $\ord_p(b)>\ord_p(n^{d-1})$ by \eqref{net}, and hence 
for any $i,$ $1\leq i\leq d-1$, the equality 
$$\ord_p\big(\GCD(b,n^i)\big)=\ord_p(n)i$$ holds.
It follows now from \eqref{brose} and \eqref{barsuk} that
$$\ord_p(b_i)=\ord_p(n)-\ord_p\big(\GCD(b,n^i)\big)=
\ord_p(n)(d-1-i)+\ord_p(n)-\ord_p(a),$$ implying that
$$ b_i\geq  p^{\ord_p(n)(d-1-i)+\ord_p(n)-\ord_p(a)}\geq  p^{\ord_p(n)(d-1-i)+1}\geq  p^{(d-1-i)+1}=
p^{d-i}.$$ 
Similarly,  since $n\nmid b$,
there exists $q\in \mathrm{rad}(n)$ such that 
$\ord_q(n)-\ord_q(b)>0$ implying by \eqref{brose} and \eqref{barsuk} that that  for any $i,$ $1\leq i\leq d-1,$  the inequality 
$ a_i\geq  q^{i}$  holds. Since $p\geq 2,$ $q\geq 2,$ this proves \eqref{neru}. 

In order to establish now the required bound, observe that
since $$A^{\circ d}=N\circ U_i\circ Y_i,$$ it follows from Corollary \ref{iter} that if \eqref{en1} or \eqref{en3} holds, then $a_i\leq 2n$. On the other hand, since $X=X_i\circ W_i$, if \eqref{en2} or \eqref{en3} holds, then $b_i\leq lb_i\leq 2n$, by Theorem \ref{ui}. 
Thus, for any $i,$ $1\leq i \leq d-1$, the inequality 
$$\min\{a_i,b_i\}\leq 2n$$ 
holds. On the other hand,  it follows from \eqref{neru} that for $i_0=\lfloor{d/2\rfloor}$ the inequality 
$$\min\{a_i,b_i\}\geq 2^{\lfloor{d/2\rfloor}}$$ holds. 
 Therefore, $2^{\lfloor{d/2\rfloor}}\leq 2n,$ implying that  
$2^{d/2}\leq 2\sqrt{2}n$. Thus, \linebreak
$d/2\leq \log_2n+3/2$ and $d\leq 2\log_2n+3.$ \qed

\vskip 0.2cm

\noindent{\it Proof of Theorem \ref{exi}.}
Observe first that if $X\in \f E(B)$ is a semiconjugacy from $B$ to $A$, then $A$ is defined in a unique way since the equalities 
$$A\circ X=X\circ B, \ \ \   \t A\circ X=X\circ B$$ imply the equality $A\circ X=\t A\circ X$ which in  turn implies the equality $A=\t A$. In particular, this implies that 
for any $X_1,X_2\in \f E(B)$ such that $X_2=\mu\circ X_1$ for some polynomial $\mu$ of degree one the corresponding polynomials  $A_1, A_2\in \f F(B)$ are conjugated.
Further, for any $A\in \f F(B)$ there exists $X$ such that
\be \l{kone} A\circ X=X\circ B\ee and $X$ is not a polynomial in $B$, since 
equalities \eqref{kone} and $X=\t X\circ B^{\circ s}$ imply the equality 
$$A\circ \t X=\t X\circ B.$$
Finally, 
 if $X_1,X_2\in \f E(B)$ and $\deg X_1=\deg X_2$, then the corresponding polynomials in $A_1, A_2\in \f F(B)$ are conjugated, since
Theorem \ref{j+} and Theorem \ref{zz} imply that there exists a polynomial $\mu$ of degree one such that $X_2=\mu\circ X_1$.

Let $X$ be an element of $\f E(B)$ and $X=\nu\circ z^{l} \circ W$ its representation  from Proposition \ref{predl}. Then it follows  from Proposition \ref{suks} that, unless $X$ is a polynomial in $B$,
the inequality $d\leq 2\log_2n+3$ holds.
Since, in addition, for the number $l$ the inequality $l<n$ holds, this implies that up to the change 
$X\rightarrow \mu\circ X$, where $\mu$ is a polynomial of degree one, there exists at most a finite number of elements of $\f E(B)$ which are not polynomials in $B.$ Applying to these  polynomials  recursively Theorem \ref{resh} we obtain polynomials $X\in \f E(B)$ and $A\in \f F(B)$ which satisfy the conclusion of the theorem.    
\qed

\br Since the degree of the polynomial of $X$ from Theorem \ref{exi} is equal to the least common multiple of degrees of all polynomials from $\f E(B)$ which are not polynomials in $B$, it follows from  Proposition  \ref{predl} and Proposition  \ref{suks} that $\deg X$ is bounded by the number $\psi(n)n^{2\log_2n+3},$ where 
$\psi(n)$ denotes the least common multiple of all numbers less than $n$ and coprime with $n$. In particular, 
$$c(n)\leq (n-1)!n^{2\log_2n+3}.$$

\er

\bc \l{zai} Let $B$ be a polynomial of degree at least two. Then  there exists at most a finite number of conjugacy classes of polynomials $A$ such that $A\leq B$.
\ec
\pr If $B$ is non-special, then the corollary follows from Theorem \ref{exi}.  For special $B$ the corollary follows Theorem \ref{sem}. \qed

\bc
Each  equivalence class of the relation $\sim$ contains at most a finite number of conjugacy classes. 
\ec
\pr Follows from Corollary \ref{zai}, since $A\sim B$ implies $A\leq B$, 

\bc [\cite{mz}]\la{mzi} Let $B$ be a non-special polynomial of degree $n\geq 2$, and $X$ and $Y$ polynomials such that $Y\circ X=B^{\circ s}$ for some $s\geq 1$. Then there exist polynomials $\t X$, $\t Y$ and $i,j\geq 0$ such that 
$$Y=B^{\circ i}\circ \t Y, \ \ \ X=\t X\circ B^{\circ j},\ \ \ {\it and} \ \ \ \t Y\circ \t X=B^{\circ \t s},$$ where $\t s$ is bounded from above by a constant which depends  
on $n$ only. 
\ec
\pr Clearly, without loss of generality we may assume that $X$ is not a polynomial in $B$. Since  $B\circ B^{\circ d}=B^{\circ d}\circ B,$ the polynomial $B^{\circ d}$ is contained in $\f E(B)$ and hence $X$ is contained in $\f E(B)$ by Corollary \ref{c1}.
Furthermore, since $\rad(\deg X)\mid \rad(n),$ it follows from 
Proposition  \ref{predl} and Proposition  \ref{suks} that there exists a polynomial $\t Y$ such that $\t Y\circ X=B^{\circ (2\log_2n+3)}.$
Therefore, if  $s> 2\log_2n+3,$ then 
$$B^{\circ s}=B^{\circ (s-2\log_2n-3)}\circ B^{\circ (2\log_2n+3)}= B^{\circ (s-2\log_2n-3)}\circ\t Y\circ X=Y\circ X,$$ implying that 
$Y=B^{\circ (s-2\log_2n-3)}\circ\t Y$. This proves the corollary, and shows that $\t s\leq 2\log_2n+3.$ \qed

\br The bound $\t s\leq 2\log_2n+3$ in Corollary \ref{mzi}  is not optimal. It was shown in \cite{mz} that in fact $\t s\leq \log_2(n+2)$ and that this last 
bound cannot be improved. For more details we refer the reader to \cite{mz}. Notice however that for applications, similar to ones given in \cite{gtz2}, the actual form of the bound for $\t s$ is not important.  
\er

\vskip 0.2cm
\noindent{\bf Acknowledgments}. The author is grateful to  the Max-Planck-Institut fuer Mathematik for the hospitality and the support.

\end{subsection}

\end{section}

\bibliographystyle{amsplain}

\end{document}